\documentclass[a4paper,12pt]{article}

\usepackage{amssymb}
\usepackage{graphicx}
\usepackage{amsmath}
\usepackage{amsfonts}

\newcommand{\singlespacing}{\let\CS=\@currsize\renewcommand{\baselinestreatch}{1.0}\tiny\CS}
\newcommand{\doublespacing}{\let\CS=\@currsize\renewcommand{\baselinestreatch}{1.5}\tiny\CS }

\newtheorem{thm}{Theorem}[section]
\newtheorem{cor}{Corollary}[section]

\newtheorem{exam}{Example}[section]
\newtheorem{lem}{Lemma}[section]

\newtheorem{remark}{Remark}[section]

\numberwithin{equation}{section}

\begin{document}
\begin{center}
\textbf{\Large {The Critical Point Equation on Kenmotsu and almost Kenmotsu manifolds}}
\end{center}
\centerline{Dhriti Sundar Patra$^1$, Amalendu Ghosh$^2$  and Arindam Bhattacharyya$^3$}

\newtheorem{Theorem}{\quad Theorem}[section]
\newtheorem{Definition}[Theorem]{\quad Definition}
\newtheorem{Corollary}[Theorem]{\quad Corollary}
\newtheorem{Lemma}[Theorem]{\quad Lemma}
\newtheorem{Example}[Theorem]{\emph{Example}}
\newtheorem{Proposition}[Theorem]{Proposition}
\numberwithin{equation}{section}
\noindent\\
\textbf{Abstract:} {In this paper, we have studied the critical point equation (shortly, CPE) within the frame-work of Kenmotsu and almost Kenmotsu manifolds. First, we prove that a complete Kenmotsu metric satisfies the CPE is Einstein and locally isometric to the hyperbolic space $\mathbb{H}^{2n+1}$.  In case of Kenmotsu manifolds, it is possible to determine the potential function explicitly (locally). We also provide some examples of Kenmotsu and almost Kenmotsu manifolds that satisfies the CPE.}

\noindent\\
\textbf{Mathematics Subject Classification 2010}: 53D15; 53C25, 53C21.

\noindent\\
\textbf{Keywords}: Total scalar curvature functional, the critical point equation, Kenmotsu manifold, almost Kenmotsu manifold, generalized nullity distribution.

\section{Introduction}
A classical problem in differential geometry is to find Riemannian metrics on a given compact
manifold that provides constant curvature. In this sense, it is interesting to study the critical points
of the total scalar curvature functional $\mathcal{S} :\mathcal{M} \longrightarrow  R$ given by
\begin{eqnarray*}
&\mathcal{S}(g) =\int_{M}r_{g}dv_{g},
\end{eqnarray*}
defined on a compact orientable Riemannian manifold $(M^n,g)$, where $\mathcal{M}$ denote the set of all Riemannian metrics on $(M^n,g)$ of unit volume, $r_{g}$ the scalar curvature and $dv_{g}$ the volume form (determined by the metric and orientation). The functional $\mathcal{S}$ restricted over $\mathcal{M}$ is known as the Einstein-Hilbert functional and its critical points are the Einstein metrics (for details, see Chapter $2$ in \cite{AB}).\\

From the classical Yamabe problem it follows that there exists many Riemannian metrics on a compact manifold that has constant scalar curvature. Let $\mathcal{C} =\{g \in M|r_{g}$ = constant$\}$. The Euler-Lagrangian equation of Hilbert-Einstein functional restricted to $\mathcal{C}$ on a given compact oriented manifold $(M, g)$ can be written as the following critical point equation (shortly, CPE)
\begin{eqnarray}\label{1.1}
&\mathcal{L}^{*}_g(\lambda)= {Ric}^{o}_{g},
\end{eqnarray}
where ${Ric}^{o}_{g}$ denotes the traceless Ricci tensor of $M$. Here $\mathcal{L}^{*}_g(\lambda)$ is the formal $L^2$-adjoint of the linearized scalar curvature operator $\mathcal{L}_g(\lambda)$ and is defined as
\begin{eqnarray}\label{1.2}
&\mathcal{L}^{*}_g(\lambda) = -(\Delta_g\lambda)g + Hess_{g} \lambda - \lambda Ric_{g},
\end{eqnarray}
 where $\Delta_g$, $Ric_{g}$, and $Hess\lambda$ are respectively the Laplacian, the Ricci tensor, and the Hessian of the smooth function $\lambda$ on $M$. The function $\lambda $ is known as the potential function. Therefore, from now on, we consider a metric $g$ with a non-trivial potential function $\lambda$ as a solution of the CPE and is denoted by $(g,\lambda)$.
Using (\ref{1.2}), one can express the equation (\ref{1.1}) in the following form
\begin{eqnarray}\label{1.3}
&Hess\lambda + (\frac{r}{n-1}g - Ric_{g})\lambda =  Ric_{g} - \frac{r}{n}g.
\end{eqnarray}
\begin{remark}
 We note that if $\lambda$ is constant in the equation (\ref{1.3}), then $\lambda = 0$ and $g$ becomes Einstein. Further, tracing (\ref{1.3}) we deduce $\Delta_{g}\lambda= - \frac{r}{n-1}\lambda$. From which it follows that $\lambda$ is an eigenfunction of the Laplacian. Since the Laplacian has non-positive spectrum, the scalar curvature must be positive.
\end{remark}

A. Besse first conjectured that the solution of the CPE must be Einstein (cf. \cite{AB}, p. 128). Since then, we find many literatures regarding the solutions of the CPE. For instance, Yun-Chang-Hwang proved \cite{YCH} that if $(g,\lambda)$ is a non-trivial solution of the CPE on an $n$-dimensional compact Riemannian manifold $M$ and satisfies one of the following conditions $(i)$ Ricci tensor is parallel $(ii)$ curvature tensor is harmonic or $(iii)$ $g$ has vanishing conformal curvature tensor, then $(M,g)$ is isometric to a standard sphere. Extending all these, Barros and Ribeiro Jr \cite{ABC} proved that the CPE conjecture is also true for half conformally flat. Further in \cite{HC}, Hwang proved that the CPE conjecture is valid under certain conditions on the bounds of the potential function $\lambda$. Recently, Nato \cite{Na} obtained a necessary and sufficient condition on the norm of the gradient of the potential function for a CPE metric to be Einstein.\\

So far, we see that the attempts have been made to confirm the CPE
conjecture satisfying either, some curvature conditions or, some
conditions on the potential functions on Riemannian manifolds.
However, the CPE has not yet been considered on some other types of
Riemannian manifolds, for instance odd-dimensional Riemannian
manifolds. Hence it desrves special attention to consider the CPE on
certain class of almost contact metric manifolds. In this direction,
Ghosh-Patra considered the $K$-contact metrics that satisfy the
critical point equation \cite{GP} and proved that the CPE conjecture
is true for this class of metrics. This intrigues us to consider the
CPE on some other almost contact metrics. Here we characterize the
solution of the CPE on certain classes of almost Kenmotsu manifolds.
The organization of the paper is as follows. After some
preliminaries on Kenmotsu and almost Kenmotsu manifolds in sect$.$
$2$, we consider the CPE on Kenmotsu manifolds in sect$.$ $3$ and we
prove that a complete Kenmotsu metric is Einstein and locally
isometric to the hyperbolic space $\mathbb{H}^{2n+1}$. Finally, we
have studied the CPE within the framework of $(\kappa,\mu)'$-almost
Kenmotsu manifold and generalized $(\kappa,\mu)$- almost Kenmotsu
manifold.

\section{Preliminaries}
In this section, we recall some basic definitions and formulas on almost Kenmotsu manifold which we shall use in the sequal. A $(2n + 1)$-dimensional smooth manifold $M$ is said to be an almost contact metric manifold if it admits a $(1,1)$ tensor field $\varphi$, a unit vector field $\xi$ (called the Reeb vector field) and a $1$-form $\eta$ such that
\begin{eqnarray}\label{2.1}
\varphi^2 X = -X + \eta(X) \xi, \hskip 0.3cm \eta(X) = g(X, \xi),
\end{eqnarray}
for any vector field $X$ on $M$. A Riemannian metric $g$ is said to be an associated (or compatible) metric if it satisfies
\begin{eqnarray}\label{2.2}
g(\varphi X, \varphi Y) = g(X,Y) - \eta(X)\eta(Y),
\end{eqnarray}
for all vector fields $X$, $Y$ on $M$. From (\ref{2.1}) and
(\ref{2.2}) it is easy to verify that $\varphi\xi = 0, $ $\eta \circ
\varphi = 0.$ An almost contact manifold $M(\varphi, \xi, \eta)$
together with a compatible metric $g$ is known as almost contact
metric manifold. On almost contact metric manifolds one can always
define a fundamental $2$-form $\Phi$ by $\Phi(X,Y) = g(X,\varphi Y)$
for all vector fields $X$, $Y$ on $M$. An almost contact metric
structure $M (\varphi, \xi, \eta, g)$ is said to be contact metric
if $\Phi = d\eta$, and is said to be almost Kenmotsu manifold if
$d\eta = 0$ and $d\Phi = 2\eta\wedge\Phi$. Further, an almost
contact metric structure is said to be normal if
$[\varphi,\varphi](X,Y) + 2d\eta(X,Y)\xi = 0$, where
$[\varphi,\varphi](X,Y) = [\varphi X,\varphi Y] + \varphi^2[X,Y]-
\varphi[\varphi X,Y] -\varphi [X,\varphi Y]$ for all vector fields
$X$, $Y$ on $M$. A normal almost Kenmotsu manifold is said to be a
Kenmotsu manifold and the normality condition is given by
$(\nabla_{X}\varphi)Y = g(\varphi X,Y)\xi - \eta(Y)\varphi X$ for
all vector fields $X$, $Y$ on $M$. On a Kenmotsu manifold \cite{K}:
\begin{eqnarray}\label{2.2A}
&\nabla_{X}\xi = X - \eta(X)\xi,\\
&R(X, Y)\xi = \eta(X)Y - \eta(Y)X,\label{2.3A}\\
&Q\xi = -2n\xi, \label{2.4A}
\end{eqnarray}
for all vector fields $X$, $Y$ on $M$, where $\nabla$ the operator
of covariant differentiation of $g$, $R$ the curvature tensor of $g$
and $Q$ the Ricci operator associated with the $(0,2)$ Ricci tensor
$Ric_g$ given by $Ric_g(X,Y)=g(QX,Y)$ for all vector fields $X$, $Y$
on $M$. In \cite{K}, it is proved that a Kenmotsu manifold
$M^{2n+1}$ is locally a warped product $I\times_{f}N^{2n}$, where
$I$ is an open interval with coordinate $t$, $f = ce^{t}$ is the
warping function for some positive constant $c$ and $N^{2n}$ is a
K\"{a}hlerian manifold. A larger class of this manifold is known as
$(0, \beta)$-Kenmotsu manifold (or simply $\beta$-Kenmotsu
manifold). An almost contact metric manifold $M$ is said to be
trans-Sasakian if there exist two functions $\alpha$ and $\beta$ on
$M$ such that
\begin{eqnarray}
    (\nabla_X\varphi)Y = \alpha(g(X, Y)\xi - \eta(Y)X) + \beta(g(\varphi X, Y)\xi - \eta(Y)\varphi X).
\end{eqnarray}
for any vector fields $X$, $Y$ on $M$. If $\alpha =0$, then $M$ is
said to be $\beta$-Kenmotsu manifold. Kenmotsu manifolds appears as
a particular case of $\beta$-Kenmotsu manifolds, for $\beta = 1$.

On an almost Kenmotsu manifold $M^{2n+1}(\varphi,\xi,\eta,g)$, we now define two operators $h$ and $l$ by $ h = \frac{1}{2}\pounds_{\xi}\varphi $ and $l= R(., \xi )\xi$ on $M$, where $R$ denotes the curvature tensor and $\pounds$ is the Lie differentiation. These two symmetric tensors of type $(1,1)$ satisfy
\begin{eqnarray}\label{2.3}
h\xi = h'\xi = 0, \hskip 0.3cm Tr~h = Tr~(h')=0, \hskip 0.3cm h\varphi = - \varphi h.
\end{eqnarray}
On an almost Kenmotsu manifold the following formula is valid
\cite{DP, DP2}:
\begin{eqnarray}\label{2.4}
\nabla_{X}\xi = - \varphi^2 X - \varphi h X,
\end{eqnarray}
for any vector field $X$ on $M$.


An almost Kenmotsu manifold $M^{2n+1}(\varphi,\xi,\eta,g)$ is said to be a generalized $(\kappa,\mu)$-almost Kenmotsu manifold if $\xi$ belongs to the generalized $(\kappa,\mu)$-nullity distribution, i.e.,
\begin{eqnarray}\label{2.10}
R(X,Y)\xi = \kappa[\eta(Y)X - \eta(X)Y] + \mu[\eta(Y)hX - \eta(X)hY],
\end{eqnarray}
for all vector fields $X$, $Y$ on $M$, where $\kappa$, $\mu$ are smooth functions on  $M$.
An almost Kenmotsu manifold $M^{2n+1}(\varphi,\xi,\eta,g)$ is said to be a generalized $(\kappa,\mu)'$-almost Kenmotsu manifold if  $\xi$ belongs to the generalized $(\kappa,\mu)'$-nullity distribution, i.e.,
\begin{eqnarray}\label{2.7}
R(X,Y)\xi = \kappa[\eta(Y)X - \eta(X)Y] + \mu[\eta(Y)h' X - \eta(X)h' Y],
\end{eqnarray}
for all vector fields $X$, $Y$ on $M$, where $\kappa$, $\mu$ are smooth functions on  $M$. Moreover, if both $\kappa$ and $\mu$ are constants in Eq$.$ (\ref{2.7}), then $M$ is called a $(\kappa,\mu)'$-almost Kenmotsu manifold.
Classifications of almost Kenmotsu manifolds with $\xi$ belong to $(\kappa$, $\mu)$-nullity distribution and $(\kappa,\mu)'$-nullity distribution have done by several authors. For more details, we refer the reader to \cite{DP,DP2,DP3,W,WL}.
On generalized $(\kappa,\mu)$ or $(\kappa,\mu)'$-almost Kenmotsu manifold with $h\neq0$, the following relations hold (see \cite{DP2})
\begin{eqnarray}\label{2.8}
h'^2 = (\kappa+1)\varphi^2~~~~~or,~~equivalently~~~~~ h^2 = (\kappa+1)\varphi^2.
\end{eqnarray}
\begin{eqnarray}\label{3.3}
Q\xi=2n\kappa \xi
\end{eqnarray}

Let $\emph{D}=ker(\eta)$ be the distribution and $X \in\emph{D}$  be an eigenvector of $h'$ with eigenvalue $\sigma$. It follows from (\ref{2.8}) that
 $\sigma ^2=-(\kappa+1)$ and therefore $\kappa\leq-1$ and $\sigma=\pm\sqrt{-\kappa-1}$. The equality holds if and only if $h=0$ (equivalently, $h'=0$). Thus, $h'\neq0$ if and only if $\kappa<-1$.

We now recall the notion of warped product manifolds for our later use. Let $(N, J, \tilde{g})$ be an almost Hermitian manifold and consider the warped product $M = \mathbb{R}\times_{f}N$ with the metric $g = g_0 + f^2\tilde{g}$, where $f$ is a positive function on $\mathbb{R}$ and $g_0$ is the standard metric on $\mathbb{R}$. define $\eta = dt$, $\xi = \frac{\partial}{\partial t}$ and the tensor field $\varphi$ is defined on $\mathbb{R}\times_{f}N$ by $\varphi X = JX$ for any vector field $X$ on $N$ and $\varphi X = 0$ if $X$ is tangent on $\mathbb{R}$. Then it is easy to testify that $M$ admits an almost contact metric structure.
An interesting characterization of an almost Kenmotsu manifold through the warped product of a real line and an almost Hermitian manifold is given by the following (see \cite{ABC})
\begin{lem}\label{lem2.1}
Let $N$ be an almost Hermitian manifold. Then the warped product $\mathbb{R}\times_{f}N$ is a $(0, \beta)$-trans Sasakian msnifold, with $\beta = \frac{f^{'}}{f}$ if and only if $N$ is K$\ddot{a}$hlerian.
\end{lem}
\section{Kenmotsu manifold satisfying the CPE}

We assume that a Riemannian manifold $(M^{2n+1}, g)$ satisfies the the CPE. Then the equation (\ref{1.1}) can be exhibited as
\begin{eqnarray}\label{3.3A}
&\nabla_{X}D\lambda=(\lambda+1)QX+(\triangle_{g}\lambda-\frac{r}{2n+1})X,
\end{eqnarray}
for any vector field $X$ on $M$. Now, tracing (\ref{1.1}) we obtain $\triangle_{g}\lambda = -\frac{r\lambda}{2n}$. Then the Eq$.$ (\ref{3.3A}) transforms into
\begin{eqnarray}\label{3.2}
\nabla_{X}D\lambda = (\lambda + 1) QX + fX,
\end{eqnarray}
for any vector field $X$ on $M$ and $f=-r(\frac{\lambda}{2n} + \frac{1}{2n+1})$.
Taking the covariant derivative of (\ref{3.2}) along an arbitrary vector field $Y$ on $M$, we get
\begin{eqnarray*}
\nabla_{Y}(\nabla_{X}D\lambda) = (Y\lambda)QX  + (\lambda
+1)\{(\nabla_{Y}Q)X + Q(\nabla_{Y}X)\} + (Yf)X + f\nabla_{Y}X,
\end{eqnarray*}
for any vector field $X$ on $M$. Applying the preceding equation and (\ref{3.2}) in the well known expression of the curvature tensor
$R(X,Y) = [\nabla_{X},\nabla_{Y}]-\nabla_{[X,Y]},$ we deduce
\begin{eqnarray}\label{3.1}
&R(X,Y)D\lambda = (X\lambda)QY - (Y\lambda)QX + (\lambda +1)\{(\nabla_{X}Q)Y -(\nabla_{Y}Q)X\} \nonumber\\
 & + (Xf)Y - (Yf)X
\end{eqnarray}
for any vector fields $X$, $Y$ on $M$.\\

First, we recall the following formula which is valid for any Kenmotsu manifold
\begin{eqnarray}\label{3.1AA}
(\nabla_{\xi}Q)Y = -2QY -4n Y,
\end{eqnarray}
for any vector field $Y$ on $M$. The complete proof of this can be found in \cite{AG}.\\

We now in position to characterize the solution of the CPE within the framwork of Kenmotsu manifold.
\begin{thm}
Let $M^{2n+1}(\varphi,\xi,\eta,g)$ be a complete Kenmotsu manifold. If $(g,\lambda)$ satisfies the CPE, then $g$ is Einstein and locally isometric to the hyperbolic space $\mathbb{H}^{2n+1}$. Further the potential function is locally given by $\lambda = A~cosh~t + B~sinh~t $, where $A$, $B$ are constants on $M$.
\end{thm}
\textbf{Proof: }Taking covariant derivative of \eqref{2.3A} over an arbitrary vector field $X$ on $M$ and using \eqref{2.2A} we have
\begin{eqnarray}\label{3.8BA}
(\nabla_{X}Q)\xi = -QX -2n X.
\end{eqnarray}
Now, substituting $\xi$ by $Y$ in \eqref{3.1} and using \eqref{2.4A}, \eqref{3.1AA}, \eqref{3.8BA}, we obtain
\begin{eqnarray}\label{3.7B}
&R(X,\xi)D\lambda = -2n(X\lambda)\xi - (\xi\lambda)QX + (\lambda+1)(QX + 2nX) \nonumber\\
&+ (Xf)\xi - (\xi f)X,
\end{eqnarray}
for any vector field $X$ on $M$. Further, from (\ref{2.3A}), we deduce $R(X,\xi)Z = g(X,Z)\xi - \eta(Z)X$. By virtue of this the foregoing equation reduces to
\begin{eqnarray}\label{3.7BB}
&\{2n(\lambda+1) - \xi f + \xi\lambda\}X - (\xi\lambda - \lambda - 1)QX \nonumber\\
&+ (Xf)\xi - (2n+1)(\xi\lambda)\xi = 0.
\end{eqnarray}
Now, taking $g$-trace of \eqref{3.1AA} yields $\xi r= -2 (r+2n(2n+1))$. Since $(g,\lambda)$ is a solution of the the CPE, the scalar curvature $r$ of $g$ is constant and hence $r=-2n(2n+1)$. Making use of this in $f=-r(\frac{\lambda}{2n} + \frac{1}{2n+1})$, we deduce $f=(2n+1)\lambda+2n$ and hence
\begin{eqnarray}\label{3.5B}
Xf = (2n+1)(X\lambda).
\end{eqnarray}
Since $\nabla_{\xi}\xi=0$ (follows from \eqref{2.4A}) and $\xi\lambda = g(\xi,D\lambda)$, taking into account \eqref{2.4A} and \eqref{3.2}, we deduce
\begin{eqnarray}\label{3.6B}
\xi(\xi\lambda)=g(\nabla_{\xi}D\lambda,\xi)=\lambda.
\end{eqnarray}
Making use of \eqref{3.5B}, the equation \eqref{3.7BB} transform into
\begin{eqnarray}\label{3.8B}
(\xi\lambda - \lambda - 1)\{QX+2nX\}=0,
\end{eqnarray}
for all vector field $X$ on $M$. If possible, let $\xi\lambda=\lambda +1$ in some open set $\mathcal{O}$ of $M$. Then from \eqref{3.6B} we have $\xi\lambda = \lambda$. These two gives $\lambda +1=\lambda$, i.e., $1=0$, which is absurd. Hence $QX=-2n X$ for all vector field $X$ on $M$ and $M$ is Einstein. Then \eqref{3.2} can be exhibited as
\begin{eqnarray}\label{3.8B}
\nabla_{X}D\lambda = \lambda X,
\end{eqnarray}
for all vector field $X$ on $M$. Applying Tashiro's theorem \cite{Tas} we conclude that $M$ is isometric to the hyperbolic space $\mathbb{H}^{2n+1}$. \\

Next, we determine the potential function $\lambda$.
It is known \cite {K} that a Kenmotsu manifold is locally a warped product $I\times_{f}N^{2n}$ , where $I$ is an open interval of the real line, $N^{2n}$ is a K$\ddot{a}$hler manifold, $f^2 = ce^{2t}$ and $t$ is the coordinate of $I$. Since $\xi = \frac{\partial}{\partial t}$, it follows from \eqref{3.8B} that
\begin{eqnarray*}
\frac{d^2\lambda}{dt^2} = \lambda.
\end{eqnarray*}
It's solution can be exhibited as $\lambda = A~cosh~t + B~sinh~t $, where $A$, $B$ are constants on $M$. This completes the proof.~~~~~~~~~~~~~~~~~~~~~~~~~~~~~~~~~~~~~~~~~~~~~~~~~~~~~~~~~~~~~~~~~~~~~$\square$
\begin{exam}
Let $(N,J, g_0)$ be a K$\ddot{a}$hler manifold of dimension $2n$. Consider the warped product $(M, g) = (\mathbb{R}\times _{\sigma} N, dt^2 + \sigma^2 g_0)$, where $t$ is the coordinate on $\mathbb{R}$.  We set $\eta = dt$, $\xi = \frac{\partial}{\partial t}$ and the tensor field $\varphi$ is defined on $\mathbb{R}\times_{\sigma}N$ by $\varphi X = JX$ for vector field $X$ on $N$ and $\varphi X = 0$ if $X$ is tangent to $\mathbb{R}$. Then the warped product $\mathbb{R}\times_{\sigma}N$, $\sigma^2 = ce^{2t}$, with the structure $(\varphi,\xi,\eta,g)$ is a Kenmotsu manifold (see \cite{K}). In particular, if we take $N = CH^{2n}$, then $N$ being Einstein, the Ricci tensor of $M$ becomes $S = -2n g$. Define $\lambda(t) = ke^t, k>0$. Then it is easy to verify that $ (g, \lambda)$ is a solution of the CPE. Similarly, we may also construct many examples that satisfies the CPE by taking different potential functions on the warped product.
\end{exam}
\section{Almost Kenmotsu manifold satisfying the CPE}
Before proving our main results we now recall some lemmas for our later use.
\begin{lem}\label{lem2.2}
(Proposition $3.2$ of \cite{DP3}) Let $M^{2n+1}(\varphi,\xi,\eta,g)$ be an almost Kenmotsu manifold with $\xi$ belonging to the generalized $(\kappa,\mu)'$-nullity distribution and $h\neq0$. Then
\begin{eqnarray}\label{A}
\xi(\lambda)=-\lambda (\mu+2), ~~~~~~~\xi(\kappa)=-2(\kappa+1)(\mu+2).
\end{eqnarray}
\end{lem}
Using Proposition $4.2$ of \cite{DP2} and Theorem $5.1$ of \cite{DP3}, Wang-Liu\cite{WL2} obtained the expression of Ricci operator on almost Kenmotsu manifold when $\xi$ belonging to the generalized $(\kappa,\mu)'$ or $(\kappa,\mu)'$-nullity distribution and $h'\neq0$.
\begin{lem}\label{lem2.3}
(Lemma $3.3$ of \cite{WL2}) Let $M^{2n+1}(\varphi,\xi,\eta,g)$ be a generalized $(\kappa,\mu)'$-almost Kenmotsu manifold with $h'\neq0$. For $n>1$, the Ricci operator $Q$ of $M$ can be expressed as
\begin{eqnarray}\label{B}
QX = -2nX + 2n(\kappa+1)\eta(X)\xi -[\mu-2(n-1)h']X,
\end{eqnarray}
for any vector field $X$ on $M$. Further, if $\kappa$ and $\mu$ are constants and $n\geq1$, then $\mu=-2$ and hence
\begin{eqnarray}\label{2.9}
QX = -2nX + 2n(\kappa+1)\eta(X)\xi - 2nh'X,
\end{eqnarray}
for any vector field $X$ on $M$. In both the cases, the scalar curvature of $M$ is $2n(\kappa-2n)$.
\end{lem}
In the similar way, using the Theorem $4.1$ of \cite{DP3}, Wang-Liu\cite{WL2} also obtained the expression of Ricci operator on almost Kenmotsu manifold when $\xi$ belonging to the generalized $(\kappa,\mu)$-nullity distribution and $h\neq0$.
\begin{lem}\label{lem2.4}
(Lemma $3.4$ of \cite{WL2}) Let $M^{2n+1}(\varphi,\xi,\eta,g)$ be a generalized $(\kappa,\mu)$-almost Kenmotsu manifold with $h'\neq0$. For $n>1$, the Ricci operator $Q$ of $M$ can be expressed as
\begin{eqnarray}\label{C}
QX = -2nX + 2n(\kappa+1)\eta(X)\xi -2(n-1)h'X + \mu hX,
\end{eqnarray}
for any vector field $X$ on $M$. Also, the scalar curvature of $M$ is $2n(\kappa-2n)$.
\end{lem}
\begin{remark}\label{rem4.1A}
We note that an Einstein almost Kenmotsu manifold satisfying any one of the hypothesis of Lemma $4.2$ or Lemma $4.3$ does not exist. Indeed, if $M$ is Einstein, then $S = \frac{r}{2n+1}g$. Since $Q\xi = 2nk\xi$ and the scalar curvature is $2n(\kappa-2n)$, we see that $k=-1$ and hence $h'= 0$. This contradicts the hypothesis $h'\neq 0$. Thus, any almost Kenmotsu manifold satisfying a nullity condition with $h'\neq0$ does not exists.
\end{remark}
\begin{lem}\label{lem4.1}
Let $M^{2n+1}(\varphi,\xi,\eta,g)$ be a generalized $(\kappa,\mu)$ or $(\kappa,\mu)'$-almost Kenmotsu manifold with $h\neq0$. If $\kappa$ is constant then we have
\begin{eqnarray}\label{D}
(\nabla_{X}Q)\xi = 2n\kappa(X-\varphi hX)-Q(X-\varphi hX),
\end{eqnarray}
for any vector field $X$ on $M$.
\end{lem}
\textbf{Proof: }
Taking covariant derivative  of (\ref{3.3}) along an arbitrary vector field $X$ on $M$ we have
$$(\nabla_{X}Q)\xi + Q(\nabla_{X}\xi)=2n\kappa \nabla_{X}\xi.$$
Making use of (\ref{2.4}) we complete the proof.~~~~~~~~~~~~~~~~~~~~~~~~~~~~~~~~~~~~~~~~~~~~~~~~~~~~~~~~~~~~~$\square$\\


Next, by using the above lemmas we prove

\begin{thm}
Let $M^{2n+1}(\varphi, \xi, \eta, g)$ be a $(\kappa,\mu)'-$almost Kenmotsu manifold with $h'\neq0$. If $(g,\lambda)$ be a non-trivial solution of the CPE then $M^3$ is locally isometric to the Riemannian product $\mathbb{H}^2(-4)\times \mathbb{R}$, and for $n>1$, $M$ is locally isometric to the warped products $\mathbb{H}^{n+1}(\alpha)\times_{f} \mathbb{R}^{n}$, $B^{n+1}(\alpha')\times_{f'}\mathbb{R}^n$, where $\mathbb{H}^{n+1}(\alpha)$ is the hyperbolic space of constant curvature $\alpha=-1-\frac{2}{n}-\frac{1}{n^2}$, $B^{n+1}(\alpha')$ is a space of constant curvature $\alpha'=-1+\frac{2}{n}-\frac{1}{n^2}$, $f=ce^{(1-\frac{1}{n})t}$ and $f'=c'e^{(1+\frac{1}{n})t}$, with $c$, $c'$ positive constants.
\end{thm}
\textbf{Proof: }
Taking scalar product of the equation (\ref{3.1}) with $\xi$ and using the value of scalar curvature (from Lemma $3.1$) and (\ref{3.3}) provides
\begin{eqnarray}\label{3.4}
&g(R(X,Y)D\lambda,\xi) = [(2n-1)\kappa+2n]\{(X\lambda)\eta(Y) - (Y\lambda)\eta(X)\} \nonumber\\
&+ (\lambda+1)\{g(Y,(\nabla_{X}Q)\xi) - g(X,(\nabla_{Y}Q)\xi)\},
\end{eqnarray}
for all vector fields $X$, $Y$ on $M$.
Now, making use of (\ref{D}) and (\ref{2.3}) in (\ref{3.4}), we immediately infer that
\begin{eqnarray}\label{E}
&g(R(X,Y)D\lambda,\xi) = [(2n-1)\kappa+2n]\{(X\lambda)\eta(Y) - (Y\lambda)\eta(X)\} \nonumber\\
&+(\lambda+1)\{g(Q\varphi hX,Y) - g(X,Q\varphi  hY)\},
\end{eqnarray}
for all vector fields $X$, $Y$ on $M$. Substituting $X$ by $\xi$ in (\ref{E})and then recalling (\ref{2.3}) yields
\begin{eqnarray}\label{F}
g(R(\xi,Y)D\lambda,\xi) = [(2n-1)\kappa+2n]\{(\xi\lambda)\eta(Y) - (Y\lambda)\},
\end{eqnarray}
for any vector field $Y$ on $M$. On the other hand, replacing $X$ by $\xi$ in the equation (\ref{2.7}) and then taking the scalar product of the resulting Eq$.$ with $D\lambda$ gives
\begin{eqnarray*}
g(R(\xi,Y)D\lambda,\xi) = \kappa g(D\lambda-(\xi\lambda)\xi,Y) + \mu g(D\lambda,h' Y),
\end{eqnarray*}
for any vector field $Y$ on $M$. Combining the last two equations we have
\begin{eqnarray*}
\mu h' D\lambda = -2n(\kappa+1)(D\lambda-(\xi\lambda)\xi).
\end{eqnarray*}
Now, operating the last equation by $\mu h'$ and then recalling the foregoing Eq$.$ provides
$\mu^2 h'^2 D\lambda = 4n^2(\kappa+1)^2(D\lambda-(\xi\lambda)\xi).$
By virtue of (\ref{2.1}) and (\ref{2.8}), the last equation transforms into
$(\kappa+1)\{4n^2(\kappa+1)+\mu^2\}(D\lambda-(\xi\lambda)\xi) = 0.$
Since $h'\neq0$,  $\kappa<-1$. Hence, the preceding equation reduces to
\begin{eqnarray}\label{3.5}
\{4n^2(\kappa+1)+\mu^2\}(D\lambda-(\xi\lambda)\xi) = 0.
\end{eqnarray}
As $\kappa$, $\mu$ are constants, we have either $4n^2(\kappa+1)+\mu^2=0$, or $4n^2(\kappa+1)+\mu^2\neq0$.\\

\noindent
\textbf{Case I: }In this case, we have $\kappa=-1-\frac{\mu^{2}}{4n^{2}}$. According to Proposition $4.1$ of Dileo and Pastore \cite{DP2}, $\mu=-2$ and therefore $\kappa=-1-\frac{1}{n^{2}}$. For $n=1$, $\kappa=\mu=-2 $ and therefore from Theorem $4.2$ of Dileo and Pastore \cite{DP2} we deduce that $M^3$ is locally isometric to the Riemannian product $\mathbb{H}^2(-4)\times \mathbb{R}$ and for $n>1$, $M$ is locally isometric to the warped products $\mathbb{H}^{n+1}(\alpha)\times_{f} \mathbb{R}^{n}$, $B^{n+1}(\alpha')\times_{f'}\mathbb{R}^n$, where $\mathbb{H}^{n+1}(\alpha)$ is the hyperbolic space of constant curvature $\alpha=-1-\frac{2}{n}-\frac{1}{n^2}$, $B^{n+1}(\alpha')$ is a space of constant curvature $\alpha'=-1+\frac{2}{n}-\frac{1}{n^2}$, $f=ce^{(1-\frac{1}{n})t}$ and $f'=c'e^{(1+\frac{1}{n})t}$, with $c$, $c'$ positive constants. \\

\noindent
\textbf{Case II: }In this case, it follows from (\ref{3.5}) that $ D\lambda = (\xi \lambda)\xi $. Taking covariant derivative of $ D\lambda = (\xi \lambda)\xi $ along an arbitrary vector field $ X $ on $M$ and using (\ref{2.1}), (\ref{2.4}), we deduce
\begin{eqnarray}\label{3.5A} \nabla_{X}D\lambda = X(\xi \lambda)\xi + (\xi \lambda)( X - \eta(X)\xi - \varphi h X).
\end{eqnarray}
Making use of (\ref{3.2}) in (\ref{3.5A}) it follows that
\begin{eqnarray*}
&(\lambda+1)QX = (\kappa-2n)(\lambda+\frac{2n}{2n+1})X + X(\xi \lambda)\xi + (\xi \lambda)( X - \eta(X)\xi - \varphi h X),
\end{eqnarray*}
for any vector field $X$ on $M$. Comparing this with (\ref{2.9}) we deduce that
\begin{eqnarray}\label{3.6}
&\{\kappa\lambda+(\xi\lambda)+(\kappa+1)\frac{2n}{2n+1}\}X  +\{2n\lambda+(\xi\lambda)+2n\}h' + X(\xi\lambda)\xi X\nonumber\\
&- \{(\xi\lambda)+2n(\kappa+1)(\lambda+1)\}\eta(X)\xi = 0,
\end{eqnarray}
for any vector field $X$ on $M$. Now, tracing (\ref{3.6}) over $X$ and noting that $Tr~h' = 0$, we have
\begin{eqnarray}\label{3.7}
&(2n+1)\{\kappa\lambda+(\xi\lambda)\}+2n(\kappa+1) + \xi(\xi\lambda) \nonumber\\
&-\{(\xi\lambda)+2n(\kappa+1)(\lambda+1)\} = 0.
\end{eqnarray}
Next, substituting $X$ by $\xi$ in the equation (\ref{3.2}) and then taking its scalar product with $\xi$ yields
$\xi(\xi\lambda) = 2n\kappa(\lambda+1)+(2n-\kappa)(\lambda+\frac{2n}{2n+1}).$
By virtue of this, equation (\ref{3.7}) takes the form
\begin{eqnarray}\label{3.8}
&\kappa\lambda+(\xi\lambda)+(\kappa+1)\frac{2n}{2n+1} = 0.
\end{eqnarray}
Therefore, operating equation (\ref{3.6}) by $\varphi$ and taking into account (\ref{3.8}) we infer that $(2n\lambda+(\xi\lambda)+2n)\varphi h' X = 0$ for any vector field $X$ on $M$. Further, the action of $\varphi$ and recalling $h' = h\circ\varphi$, the last Eq$.$ provides
\begin{eqnarray}\label{3.9}
(2n\lambda+(\xi\lambda)+2n)h'X = 0,
\end{eqnarray}
for any vector field $X$ on $M$. By virtue of (\ref{3.8}), the equation (\ref{3.9}) reduces to $(2n-\kappa)(\lambda+\frac{2n}{2n+1})h'X = 0$ for any vector field $X$ on $M$.
Since $h'$ is non-vanishing, $\kappa<-1$, it follows that $\lambda = -\frac{2n}{2n+1}$ and which is constant. Therefore, from \textbf{Remark 1.1}, it follows that $\lambda=0$ and therefore $2n=0$, a contradiction. This completes the proof.~~~~~~~~~~~~~~~~~~~~~~~~~~~~~~~~~~~~~~~~~~~~~~~~~~~~~~~~~~~~~~~~~~~~~~~~~~~~~~~~$\square$\\

Using the proof of the last theorem, we can extend the last result for a generalized $(\kappa,\mu)'-$almost Kenmotsu manifold.

\begin{cor}
Let $M^{2n+1}(\varphi, \xi, \eta, g)$, $n>1$, be a generalized $(\kappa,\mu)'-$almost Kenmotsu manifold with $h'\neq0$. If $(g,\lambda)$ be a non-trivial solution of the CPE then $g$ is locally isometric to the warped products $\mathbb{H}^{n+1}(\alpha)\times_{f} \mathbb{R}^{n}$, $B^{n+1}(\alpha)\times_{f'}\mathbb{R}^n$, where $f=ce^{(1-\frac{1}{n})t}$ and $f'=c'e^{(1+\frac{1}{n})t}$, with $c$, $c'$ positive constants.
\end{cor}
\textbf{Proof:} From Lemma \ref{lem2.3}, the scalar curvature of $M$ is $2n(\kappa-2n)$. But we know that the Riemannian manifold satisfying the CPE has constant scalar curvature (see Hwang \cite{HC}) and hence $\kappa$ is also constant. Further, we note that the Lemma \ref{lem4.1} is also valid here . Therefore, from Lemma \ref{lem2.2} it follows that $(\kappa+1)(\mu+2)=0$ . Since $\kappa<-1$, we must have $\mu=-2$. Thus, we have proved that the generalized $(\kappa,\mu)'-$almost Kenmotsu manifold reduces to $(\kappa,\mu)'-$almost Kenmotsu manifold. Hence the proof follows from the last theorem.~~~~~~~~~~~~~~~~~~~~~~~~~~~~~~~~~~~~~~~~~~~~~~~~~~~~~~~~~~~~~~~~~~~~~~~~~~~~~~~~~~~~~~~$\square$

\begin{thm}
If $M^{2n+1}(\varphi, \xi, \eta, g)$, $n>1$, be a generalized $(\kappa,\mu)-$almost Kenmotsu manifold with $h\neq0$ then the CPE has no solution on $M$.
\end{thm}
\textbf{Proof:} First, we note that the scalar curvature of $M$ is $2n(\kappa-2n)$ (follows from Lemma \ref{lem2.4}). Since the scalar curvature is constant (see Hwang \cite{HC}), $\kappa$ is also constant. Therefore, the Lemma \ref{lem4.1} is also valid here. Further, the equations (\ref{E}), (\ref{F}) are are true in this case.
Thus, putting $X=\varphi X$ and $Y=\varphi Y$ in (\ref{E}) and noting that $ g(R(\varphi X,\varphi Y)D\lambda,\xi) = 0$ (follows from (\ref{2.10})) and $h\varphi=-\varphi h$ we have
\begin{eqnarray*}
(\lambda+1)\{g(Qh\varphi^2X, \varphi Y)-g(\varphi X, Qh\varphi^2Y)\}=0,
\end{eqnarray*}
for all vector fields $X$, $Y$ on $M$. Making use of (\ref{C}) and $h\varphi=-\varphi h$ we have $(\lambda+1)\mu h^2 \varphi^2 X=0$ for any vector field $X$ on $M$. Using (\ref{2.1}), the last equation gives $(\kappa+1)(\lambda+1)\mu h^2 X=0$ for any vector field $X$ on $M$. Since $h\neq0$, $\kappa<-1$, it follows that $(\lambda+1)\mu=0$.\\

If $\lambda+1=0$, then from \textbf{Remark 1.1}, it follows that $\lambda=0$ and this leads to a contradiction. So, we suppose that $\lambda+1\neq0$ in some open set $\mathcal{O}$ of $M$. Then on $\mathcal{O}$, $\mu=0$. All our next discussion will be on $\mathcal{O}$.
Now, replacing $X$ by $\xi$ in the equation (\ref{2.10}) and then taking the scalar product of the resulting Eq$.$ with $D\lambda$ gives
\begin{eqnarray*}
g(R(\xi,Y)D\lambda,\xi) = \kappa g(D\lambda-(\xi\lambda)\xi,Y) + \mu g(D\lambda,h Y).
\end{eqnarray*}
Moreover, the result (\ref{F}) valid here also and therefore combining foregoing Eq$.$ with (\ref{F}) we have
\begin{eqnarray*}
2n(\kappa+1)(D\lambda-(\xi\lambda)\xi)+\mu h D\lambda =0.
\end{eqnarray*}
In this case, last Eq$.$ gives $(\kappa+1)(D\lambda - (\xi \lambda)\xi)=0 $. Since $h\neq0$, i.e., $\kappa<-1$, the last equation gives $D\lambda - (\xi \lambda)\xi=0$. Taking covariant derivative of $ D\lambda = (\xi \lambda)\xi $ along an arbitrary vector field $ X $ and using (\ref{2.1}), (\ref{2.4}), we deduce
$\nabla_{X}D\lambda = X(\xi \lambda)\xi + (\xi \lambda)( X - \eta(X)\xi - \varphi h X).$
By virtue of this, Eq$.$ (\ref{3.2}) reduces to
\begin{eqnarray*}
&(\lambda+1)QX = (\kappa-2n)(\lambda+\frac{2n}{2n+1})X + X(\xi \lambda)\xi + (\xi \lambda)( X - \eta(X)\xi - \varphi h X).
\end{eqnarray*}
By the virtue (\ref{C}), the last equation can written as
\begin{eqnarray}\label{3.12}
&\{\kappa\lambda+(\xi\lambda)+(\kappa+1)\frac{2n}{2n+1}\}X  +\{2(n-1)(\lambda+1) +(\xi\lambda)\}h' X + X(\xi\lambda)\xi \nonumber\\
& - \{2n(\kappa+1)(\lambda+1)+(\xi\lambda)\}\eta(X)\xi  = 0.
\end{eqnarray}
Now, tracing (\ref{3.12}) over $X$ and noting that $Tr~h' = 0$, we have
\begin{eqnarray}\label{3.13}
&(2n+1)\{\kappa\lambda+(\xi\lambda)\} -\{2n(\kappa+1)(\lambda+1)+(\xi\lambda)\}\nonumber\\
&+2n(\kappa+1) + \xi(\xi\lambda) = 0.
\end{eqnarray}
Next, substituting $X$ by $\xi$ in (\ref{3.2}) and using $Q\xi=2n\kappa \xi$ and then taking scalar product of the resulting equation with $\xi$, we get
$\xi(\xi\lambda) = 2n\kappa(\lambda+1)+(2n-\kappa)(\lambda+\frac{2n}{2n+1}).$
Making use of this, (\ref{3.13}) reduces to
\begin{eqnarray}\label{3.14}
&\kappa\lambda+(\xi\lambda)+(\kappa+1)\frac{2n}{2n+1} = 0.
\end{eqnarray}
Now, operating (\ref{3.12}) by $\varphi$ and using (\ref{3.14}) we get $\{2(n-1)(\lambda+1) +(\xi\lambda)\}\varphi h' X = 0$. Further, recalling $h' = h\circ\varphi$, $h\varphi=-\varphi h$, $h\xi=0$ and (\ref{2.1}), the last Eq$.$ gives
$\{2(n-1)(\lambda+1) +(\xi\lambda)\}hX = 0,$. By vietue of (\ref{3.14}), the last Eq$.$ transform into $\{(2n+1)(2n-\kappa-2)\lambda+2n(2n-\kappa-2)-2\} hX = 0$.
Since $h\neq0$, the foregoing Eq$.$ shows that $\{(2n+1)(2n-\kappa-2)\lambda+2n(2n-\kappa-2)-2\}=0$. This shows that $\lambda$ is constant. Thus, the \textbf{Remark 1.1} provides $\lambda=0$ and therefore $\kappa=2n-\frac{1}{n}-2\geq 1$, which is not possible by hypothesis. This competes the proof.~~~~~~~~~~~~~~~~~~~~~~~~~~~~~~~~~~~~~~~~~~~~~~~~~~~~~~~~~~~~~~~~~~~~~~~~~~~~~~~~$\square$\\

We shall now construct some examples of almost Kenmotsu manifolds that satisfies the CPE.
\begin{exam}\label{4.1exam}
Let $(N,J,\tilde{g})$ be a strictly almost K$\ddot{a}$hler Einstein manifold. We set $\eta = dt$, $\xi = \frac{\partial}{\partial t}$ and the tensor field $\varphi$ is defined on $\mathbb{R}\times_{f}N$ by $\varphi X = JX$ for vector field $X$ on $N$ and $\varphi X = 0$ if $X$ is tangent to $\mathbb{R}$. Consider a metric $g = g_{0}+ f^2 \tilde{g}$, where $f^2 = ce^{2t}$, $g_0$ is the Euclidean metric on $\mathbb{R}$ and $c$ is a positive constant. Then it is easy to verify (see \cite{DP}) that the warped product $\mathbb{R}\times_{f}N$, $f^2 = ce^{2t}$, with the structure $(\varphi,\xi,\eta,g)$ is an almost Kenmotsu manifold. Now, if we take  $\lambda = ce^t + \frac{1}{2n+1}$ on $M$, then it is easy to see that $\lambda$ is a solution of the CPE.
\end{exam}

\begin{remark}\label{rem4.1}
One can construct strictly almost K$\ddot{a}$hler structure on some product manifold. In fact, Oguro and Sekigawa (see \cite{OS}) constructed a strictly almost K$\ddot{a}$hler structure on the Riemannian product $\mathbb{H}^3\times \mathbb{R}$. By virtue of this it is possible to obtain a $5$-dimensional strictly almost Kenmotsu manifold on the warped product $\mathbb{R}\times_ {f^2} (\mathbb{H}^3 \times \mathbb{R})$, where $f^2 = ce^{2t}$.
\end{remark}

\begin{exam}\label{4.3exam}
Suppose that $(N^{2n}, J, ~\tilde{g})$ is a K$\ddot{a}$hler Einstein manifold with  $~\tilde{S} = -2n~\tilde{g}$. Let $f : \mathbb{R}\rightarrow \mathbb{R}^+$ be a smooth function defined by
\begin{eqnarray*}
    f(t) = A~sinh~t + B~cosh~t,
\end{eqnarray*}
where $A, B$ are constants, not simultaneously zero. Consider the warped product $M = \mathbb{R}\times _fN$ of dimension $2n+1$ endowed with the metric
\begin{eqnarray*}
g = dt^2 + f^2~\tilde{g}.
\end{eqnarray*}
Then it follows from \cite{PRRS} (e.g. see Lemma $1.1$) that $(M^{2n+1}, g)$ is an Einstein manifold. Defining $\xi$, $\eta$ and $\varphi$ as in the Example \ref{4.1exam} it is easy to see that $(M, g)$ admits an almost contact metric structure $(\varphi, \xi, \eta, g)$. Further, from Lemma \ref{lem2.1} it is obvious that the manifold $M$ with the almost contact structure is almost Kenmotsu (in particular, $\beta$-Kenmotsu). Now if we take $\lambda = A~cosh~t + B~sinh~t + \frac{1}{2n+1}$, then it is straightforward to testify that $(g, \lambda)$ is a solution of the CPE.
\end{exam}
\begin{remark}
From these examples, apparently the examples contradict the theorems that we have proved in this section. However, this is not true. Because the warped product manifolds that we have considered in the examples are Einstein while from Remark \ref{rem4.1A} it is obvious that any Einstein almost Kenmotsu manifold satisfying any nullity condition (that we have considered in the Theorems) with $h'\neq0$ does not exists.
\end{remark}
\noindent \textbf{Acknowledgments:} The authors are very much
thankful to the anonymous referee for valuable comments. The author
D. S. Patra is financially supported by the Council of Scientific
and Industrial Research, India (grant no. 17-06/2012(i)EU-V).

$^1$
Department of Mathematics, \\
Jadavpur University,  \\
188. Raja S. C. Mullick Road,\\
Kolkata:700 032, INDIA \\
E-mail: dhritimath@gmail.com\\

$^2$
Department of Mathematics, \\
Chandernagore College\\
Hooghly: 712 136 (W.B.), INDIA\\
E-mail: aghosh\_70@yahoo.com\\

$^3$
Department of Mathematics, \\
Jadavpur University,  \\
188. Raja S. C. Mullick Road,\\
Kolkata:700 032, INDIA \\
E-mail: bhattachar1968@yahoo.co.in\\


\begin{thebibliography}{99}
\bibitem{AB} A. Besse, \emph{Einstein manifolds}, Springer-Verlag, New York, 2008.
\bibitem{ABC} A. Barros, E. Ribeiro Jr., \emph{Critical Point Equation on four-dimensional compact manifolds}, Math. Nachr. 287(2014), 1618-1623.
\bibitem{Blair} D. E. Blair, \emph{Riemannian geometry of contact and symplectic manifolds}, Birkhauser, Boston, 2002.
\bibitem{BKP} D. E. Blair, T. Koufogiorgos, B. J. Papantoniou, \emph{Contact metric manifolds satisfying a nullity condition}, Israel Journal of Mathematics, 91(1995), 189-214.
\bibitem{DP} G. Dileo, A. M. Pastore, \emph{Almost Kenmotsu manifolds and local symmetry}, Bull. Belg. Math. Soc. Simon Stevin 14(2007), 343-354.
\bibitem{DP2} G. Dileo, A. M. Pastore, \emph{Almost Kenmotsu manifolds and nullity distributions}, J. Geom. 93(2009), 46-61.
\bibitem{DP3} A. M. Pastore, V. Saltarelli, \emph{Generalized nullity distributions on almost Kenmotsu manifolds}, International Electronic Journal of Geometry, 4(2)(2011), 168-183.
\bibitem{GP} A. Ghosh, D. S. Patra, \emph{The critical point equation and contact geometry}, J. Geom. 108(2017), 185-194.
\bibitem{AG} A. Ghosh, \emph{Ricci soliton and Ricci almost soliton within the framework of Kenmotsu manifold}, submitted.
\bibitem{OZ} Z. Olszak, \emph{On almost cosymplectic manifolds}, Kodai Math. J. 4(1981), 239-250.
\bibitem{CLM} D. Chinea, M. de Le´on, J. C. Marrero, \emph{Topology of cosymplectic manifolds}, J. Math. Pures Appl. 72(1993), 567-591 .
\bibitem{OS}T. Oguro and K. Sekigawa, \emph{Almost K¨ahler structures on the Riemannian product of a 3-dimensional hyperbolic space and a real line}, Tsukuba J. Math. 20, No.1 (1996), 151-161.
\bibitem{K} K. Kenmotsu, \emph{A class of almost contact Riemannian manifolds}, T\^{o}hoku Math. J. 24(1972), 93-103.
\bibitem{HC} S. Hwang, \emph{Critical points of the total scalar curvature functionals on the space of metrics of constant scalar curvature}, Manuscripta Math. 103(2000), 135-142.
\bibitem{Na} B. L. Nato, \emph{A note on critical point metrics of the total scalar curvature functional}, J. Math. Anal. Appl. 424(2015), 1544-1548.
\bibitem{Tas} Y. Tashiro, \emph{Complete Riemannian manifolds and some vector fields}, Trans. Amer. Math.  Soc. 117 (1965), 251-275.
\bibitem{W} Y. Wang, \emph{A Generalization of the Goldberg Conjecture for CoK¨ahler Manifolds}, Mediterr. J. Math. 13(5)(2016), 2679-2690.
\bibitem{WL} Y. Wang and X. Liu, \emph{On a type of almost Kenmotsu manifolds with harmonic curvature tensors}, Bull. Belg. Math. Soc. Simon Stevin 22(2015), 15-24.
\bibitem{WL2} Y. Wang and X. Liu, \emph{On almost Kenmotsu manifolds satisfying some nillity distributions}, Proc. Natl. Acad. Sci., India, Sect. A Phys. Sci.  86(3)(2016), 347-353.
\bibitem{PRRS} S. Pigola, M. Rigoli, M. Rimold, A. Setti, \emph{Ricci almost solitons}, Ann. Sc. Norm. Super. Pisa Cl. Sci. (5) 10(2011), 757-799.
\bibitem{YCH} G. Yun, S. Hwang, J. Chang, \emph{Total scalar curvature and harmonic curvature}, Taiwanese Journal of Mathematics,  18(5)(2014), 1439-1458.

\end{thebibliography}
\end{document}